\documentclass[11pt]{amsart}
\usepackage{geometry}                % See geometry.pdf to learn the layout options. There are lots.
\usepackage{amssymb}
\usepackage{epstopdf}
\usepackage{fullpage}
\usepackage{mathrsfs}

\title{Bounds on the number of cells \\ and the dimension of the Dressian}
\author{Rudi Pendavingh}
\date{\today}                                           % Activate to display a given date or no date

\begin{document}
\newcommand{\Dr}{\mathscr{D}}
\newcommand{\dr}{\mathcal{D}}
\newcommand{\R}{\mathbb{R}}
\newcommand{\N}{\mathbb{N}}
\newcommand{\Z}{\mathbb{Z}}
\newcommand{\Q}{\mathbb{Q}}
\newtheorem{theorem}{Theorem}
\newtheorem{lemma}{Lemma}
\newtheorem{corollary}{Corollary}
\newtheorem{conjecture}{Conjecture}
\newtheorem{definition}{Definition}
\newtheorem{problem}{Problem}
\newtheorem*{remark*}{Remark}
\newcommand{\ignore}[1]{}
\newcommand{\Rinf}{\R_{\infty}}

\begin{abstract}
The {\em Dressian} of a matroid $M$ is the set of all valuations of $M$. This Dressian $\Dr(M)$ is the support of a polyhedral complex $\dr(M)$ whose open cells correspond 1-1 with matroid subdivisions of the matroid polytope of $M$. 
We present upper bounds on the number of cells and the dimension of $\dr(M)$. For matroids $M$ of rank $r\geq 3$ on  $n$ elements we show that $$ \ln\#\dr(M)\leq  {\binom{n}{r}} O\left(\frac{\ln(n)^2}{n}\right)\text{ as }n\rightarrow\infty,\qquad\text{and}\qquad\dim \dr(M)\leq  {\binom{n}{r}}\frac{3}{n-r+3},$$ 
as well as some more detailed bounds that incorporate structural properties of such $M$. For uniform matroids $M=U(r,n)$, these upper bounds are comparable to lower bounds derived from valuations that are constructed from sparse paving matroids.
 \end{abstract}

\maketitle

\section{Introduction}
\subsection{Valuated matroids and the Dressian} A valuation of a matroid $M$ with bases $\mathcal{B}$ is a function $\nu:\mathcal{B}\rightarrow \R$ so that 
\begin{itemize}
\item[\bf (V)] for all $B, B'\in \mathcal{B}$ and all $e\in B\setminus B'$ there exists an $f\in B'\setminus B$ such that
$$\nu(B)+\nu(B')\geq \nu(B-e+f)+\nu(B'+e-f).$$
\end{itemize}
Valuated matroids originate from the work of Dress and Wenzel \cite{DressWenzel1992a} as a species of {\em matroids with coefficients} \cite{Dress1986}, and as such are a close relative of oriented matroids. Within the framework of matroids with coefficients, axiom {\bf(V)} is analogous to the chirotope axioms of oriented matroids and the symmetric base exchange axioms of matroids.

In tropical geometry, a matroid valuation is interpreted as the {\em tropical Pl\"ucker coordinates} of a {\em tropical linear space} \cite{Speyer2008}. Herrmann, Jensen, Joswig, and Sturmfels consider the collection of all valuations of a given (uniform) matroid $M$ in \cite{HJJS2009}, and name this set the {\em Dressian}  $\Dr(M)$.
As each valuation of $M$ is a function that assigns a real value to each basis of $M$, a valuation is a vector $\nu\in \R^{\mathcal{B}}$ and the Dressian of $M$ is a subset of $\R^{\mathcal{B}}$, where $\mathcal{B}$ denotes the set of bases of $M$. As it turns out, the Dressian $\Dr(M)$ is the support of a {\em fan} $\dr(M)$: a polyhedral complex whose cells are cones. To fully describe a Dressian, it suffices to enumerate its cells.

The Dressian of matroids of rank 2 is well-understood \cite{DressTerhalle1998,SpeyerSturmfels2004}: each matroid valuation in rank 2 is induced by a {\em metric tree}, and each cell of the Dressian corresponds to a tree topology.  
In \cite{HJJS2009}, detailed description of the cells of the Dressian of $U(3,n)$ for $n\leq 7$ is given. A key insight from this work is that each valuation $\nu$ of a matroid $M=U(3,n)$ is uniquely determined by the {\em tree arrangement} associated with the rank-2 valuations $\nu/e$ that arise by contracting each element $e$ of the ground set of $M$. This concise description is key to obtaining the enumeration of the cells of the Dressian computationally.

In many ways the present paper is a continuation of \cite{HJJS2009}. We take the same polyhedral perspective, and study the cells of the Dressian, but rather than an exact enumeration we aim for (asymptotic) upper bounds on the dimension and the number of cells. 
Where they lift exact results in rank 2 to rank 3, we argue that upper bounds in any fixed low rank $t$ induce upper bounds in any higher rank $r$. 

For a matroid $M$ of rank $r$ on $n$ elements, we obtain that 
$$\frac{\dim \dr(M)}{\binom{n}{r}}\leq \frac{\max\{\dim\dr(M/S), ~S\in \binom{E}{r-t}\text{ independent}\}}{\binom{n-r+t}{t}}.$$
The best bound on $\dim\dr(M)$ in general rank $r$ then follows a bound on the dimension in rank $t=3$ that is due to Speyer. The methods are inspired by asymptotic bounds on number of matroids due to Bansal, van der Pol, and the present author \cite{BPvdP2014}.
In that paper, we show that
$$\frac{\log (1+m(r,n))}{\binom{n}{r}}\leq \frac{\log (1+m(t,n-r+t))}{\binom{n-r+t}{t}}$$
where $m(r,n)$ denotes the number of matroids of rank $r$ on ground set $E=\{1,\ldots, n\}$. The argument hinges on the basic fact that contracting an independent set with $r-t$ elements in a matroid or rank $r$ gives a matroid of rank $t$. This key construct, contraction, generalizes to valuated matroids: if $\nu$ is a valuation of a matroid $M$, then contraction of an independent set $S$ of $M$ gives a valuation $\nu/S$ of $M/S$, and the map $\nu\mapsto \nu/S$ is essentially a coordinate projection for each $S$.
Where \cite{BPvdP2014} used Shearers Entropy Lemma to combine bounds in low rank $t$ efficiently to a bound in rank $r$, the present paper uses a structurally similar statement,  Lemma \ref{lem:notshearer},  which bounds the dimension of a linear subspace $L$  in terms of the dimensions of certain coordinate projections of $L$. 

For the number of cells of the Dressian, we find that 
$$\frac{\ln \#\dr(M)}{\binom{n}{r}}\leq \frac{\max\{\ln \#\dr(M/S): ~S\in \binom{E}{r-t}\text{ independent }\}}{\binom{n-r+t}{t}}O(\ln(n)).$$
Bounds on  $\#\dr(M)$ then follow the bounds in rank 2 that are immediate from the description of valuations in rank 2 in terms of metric trees. Our methods may be compared to the ones used in \cite{BPvdP2015} to bound on the number of matroids. In that paper, we showed that each matroid is uniquely determined by a stable set of the Johnson graph together with a modest amount of further information, so that on a logarithmic scale the number of matroids is dominated by the number of such stable sets. We then bounded the number of stable sets of the Johnson graph using the {\em container method}. In the present paper, we argue that each cell of the Dressian is uniquely determined by its linear hull. Then, we develop a rough analogue of the container method for collections of linear subspaces, Lemma \ref{lem:container}, that allows us to translate upper bounds on the number of cells of the Dressian in rank $t$ to upper bounds in rank $r$. In qualitative terms, Lemma \ref{lem:container} states that if certain coordinate projections of a collection $\mathcal{L}$ of linear subspaces  can be captured in few subspaces of small dimension, then all of $\mathcal{L}$ can be captured in few subspaces of small dimension. We then complete the argument by bounding the number of cells of the Dressian contained in a space of low dimension.

Finally, we work out lower bounds for uniform matroids that follow a construction of valuations from matroids. Essentially each sparse paving matroid $N$ determines a unique cell of the Dressian whose dimension equals the number of circuit-hyperplanes of $N$. This allows us to compare upper and lower bounds for uniform matroids $M=U(r,n)$. We find that
$$\frac{1}{n}\leq \frac{\dim \dr(M)}{\binom{n}{r}}\leq \frac{3}{n-r+1} \qquad\text{ and }\qquad\frac{1}{n}\leq \frac{\log_2 \#\dr(M)}{\binom{n}{r}}\leq c\frac{\ln(n)^2}{n}$$
for some sufficiently high constant $c$. We conjecture that the valuations obtained from sparse paving matroids have maximal dimension and dominate the count, so that the lower bounds will be closer to the truth.

\section{Preliminaries}
\subsection{Matroids} A {\em matroid} is a pair $M=(E,\mathcal{B})$ where $E$ is a finite set, the {\em ground set}, and $\mathcal{B}$ is a nonempty set of subsets of $E$, the {\em bases}, such that
\begin{itemize}
\item[\bf (B)] for all $B, B'\in \mathcal{B}$ and all $e\in B\setminus B'$ there exists an $f\in B'\setminus B$ such that 
 $$B-e+f\in\mathcal{B}\text{ and }B'+e-f\in \mathcal{B}.$$
\end{itemize}
It is immediate from {\bf (B)} that all $B\in \mathcal{B}$ have the same cardinality. This common cardinality is the {\em rank} of the matroid $M=(E,\mathcal{B})$.  

\subsection{Valuated matroids} 
Let $\Rinf:=\R\cup\{\infty\}$. A {\em valuated matroid} of rank $r$ on $E$  is a function $\nu: \binom{E}{r}\rightarrow \Rinf$ such that 
\begin{itemize}
\item[\bf (V)] for all $B, B'\in \binom{E}{r}$ and all $e\in B\setminus B'$ there exists an $f\in B'\setminus B$ such that
$$\nu(B)+\nu(B')\geq \nu(B-e+f)+\nu(B'+e-f).$$
\end{itemize}
The {\em finite support} of a function $\nu: \binom{E}{r}\rightarrow \Rinf$ is $$\underline{\nu}:=\left\{B\in\binom{E}{r}: \nu(B)\neq \infty\right\}.$$
It is straightforward that if $\nu: \binom{E}{r}\rightarrow \Rinf$ satisfies {\bf (V)}, then  $\underline{\nu}$ satisfies {\bf (B)}.

Valuated matroids were first defined by Dress and Wenzel \cite{DressWenzel1992a}, who gave the following characterization. 
\begin{theorem}\label{thm:dw}Let $\nu: \binom{E}{r}\rightarrow \Rinf$.
Then $\nu$ is a valuated matroid if and only if 
\begin{enumerate}
\item $\underline{\nu}$ satisfies the base exchange axiom {\bf (B)}; and 
\item  for each $S\in \binom{E}{r-2}$ and distinct elements $a,b,c,d\in E\setminus S$, the minimum of 
$$\nu(Sab)+\nu(Scd), ~ \nu(Sac)+\nu(Sbd), ~ \nu(Sad)+\nu(Sbc)$$
is attained at least two times.
\end{enumerate}
\end{theorem}
\ignore{\proof Necessity is straightforward. We prove sufficiency. Let $\nu, B, B', e$ be a counterexample with $|B\setminus B'|$ as small as possible. So conditions (1) and (2) hold for $\nu$,  and we have
$$\nu(B)+\nu(B')< \nu(B-e+f)+\nu(B'+e-f)$$
for all $f\in B'\setminus B$. It follows that  $|B\setminus B'|>1$. If $|B\setminus B'|=2$, then $$\nu(B)+\nu(B')<\min\{\nu(B-e+f)+\nu(B'+e-f), \nu(B-e+f')+\nu(B'+e-f')\}$$ where $B'\setminus B=\{f,f'\}$. Taking $S=B\cap B'$, $\{a,b,c,d\}=B\cup B'-S$, this violates assumption (2).
So $|B\setminus B'|>2$. 

Let $w:\R^E$ be such that $w[B]+w[B']=\nu(B)+\nu(B')$, and
$$w[B-e+f]-w[B'+e-f]=\nu(B-e+f)-\nu(B'+e-f)$$
for each $f\in B'\setminus B$. Such a $w$ exists, since these $|B'\setminus B|+1$ linear equations are linearly independent. By Lemma \ref{lem:veq}, the function $\nu': \binom{E}{r}\rightarrow \Rinf$ determined by $\nu':B''\mapsto \nu(B'')-w[B'']$ for all  $B''\in \binom{E}{r}$ is a valuated matroid. Moreover, we have $\nu'(B)+\nu'(B')=0$ and 
$\nu'(B-e+f)=\nu'(B'+e-f)$
for all $f\in B'\setminus B$. Also, we have $$\nu'(B-e+f)+\nu'(B+e-f)-\nu'(B)-\nu'(B'))=\nu(B-e+f)+\nu(B+e-f)-\nu(B)-\nu(B')>0$$
for all $f\in B'\setminus B$. Hence $$0<\min\{\nu'(B-e+f),\nu'(B'+e-f)\}$$ for each $f\in B'\setminus B$. In what follows, we replace $\nu$ with $\nu'$.

Let $X$ attain the minimum of 
$$\min\left\{\nu(X): ~X=B-e'+f', ~e'\in B\setminus B'-e, ~f'\in B'\setminus B \right\}$$
Since $\nu(B)+\nu(B')=0$, we have $B,B'\in \underline{\nu}$, and by assumption (1) {\bf (B)} holds for $\underline{\nu}$. Hence for each $e'\in B\setminus B'-e$ there exists an  $f'\in B'\setminus B$ so that $B-e'+f'\in \underline{\nu}$, i.e. so that $\nu(B-e'+f')<\infty$.  It follows that $\nu(X)<\infty$.
As $|X\setminus B'|<|B\setminus B'|$, we have $$\nu(X)+\nu(B')\geq \nu(X-e+f)+\nu(B'+e-f)>\nu(X-e+f)$$ for some $f\in B'\setminus X\subseteq B'\setminus B$.
Put $Y:=X-e+f$. Then $|B\setminus Y|=2<|B\setminus B'|$, hence $$\nu(B)+\nu(Y)\geq \nu(B-e+g)+\nu(Y+e-g)>\nu(Y+e-g)$$ for some $g\in Y\setminus B\subseteq B'\setminus B$. Hence for
$X'=Y+e-g$, we have $$\nu(X)=\nu(X)+\nu(B)+\nu(B')> \nu(B)+\nu(Y)>\nu(X'),$$ which contradicts the choice of $X$.
\endproof
}
Note that in the above statement, and in the remainder of the paper, we abbreviate $Sab:=S\cup\{a,b\}$ etc. to avoid clutter. 
\subsection{Dressians}
Given a finite set $E$ and a natural number $r\leq |E|$, we denote
$$\Dr(r,E):=\left\{\nu: \nu\text{ is a valuated matroid of rank $r$ on $E$}\right\}$$
We will consider $\nu, \nu'\in\Dr(r,E)$ as {\em combinatorially equivalent}, notation $\nu\sim\nu'$, if $\nu$ and $\nu'$ satisfy the conditions of Theorem \ref{thm:dw} in the same way; that is, if 
\begin{enumerate}
\item $\underline{\nu}=\underline{\nu'}$, and 
\item for each $S\in \binom{E}{r-2}$ and distinct elements $a,b,c,d\in E\setminus S$, we have $$\nu(Sac)+\nu(Sbd)=\nu(Sad)+\nu(Sbc) ~ \Longleftrightarrow ~ \nu'(Sac)+\nu'(Sbd)=\nu'(Sad)+\nu'(Sbc)$$
\end{enumerate}
We will write $$Z(r,E):=\left\{(S,ab,cd): S\in \binom{E}{r-2}, a,b,c,d\in E\setminus S\text{ distinct }\right\}$$ for the set of locations of condition (2) and 
$$[\nu]:=\left\{(S,ab,cd)\in Z(r,E): \nu(Sac)+\nu(Sbd)=\nu(Sad)+\nu(Sbc)\right\}$$
for each $\nu\in \Dr(r,E)$. With this notation,  the definition of combinatorial equivalence can be stated as 
$$\nu\sim \nu' \qquad :\Longleftrightarrow \qquad \underline{\nu}=\underline{\nu'} \text{ and }[\nu]=[\nu']$$

If $\mathcal{B}\subseteq \binom{E}{r}$, then any function $\nu: \mathcal{B}\rightarrow \R$ extends uniquely to a function $\overline{\nu}:\binom{E}{r}\rightarrow \Rinf$ whose finite support is $\mathcal{B}$, by setting
$$\overline{\nu}(B)=\left\{\begin{array}{ll} \nu(B)&\text{if }B\in \mathcal{B}\\ \infty &\text{otherwise}\end{array}\right.$$
If $M=(E, \mathcal{B})$ is a matroid, then $\nu:\mathcal{B}\rightarrow \R$ is a {\em valuation of $M$} if $\overline{\nu}$ is a valuated matroid.
The {\em Dressian} of $M$ is 
$$\Dr(M):=\{\nu\in \R^{\mathcal{B}}: \nu \text{ is a valuation of }M\}.$$

Two valuations  $\nu,\nu'$ are combinatorially equivalent if their extensions are, that is, $\nu\sim\nu'$ if and only if $\overline{\nu}\sim\overline{\nu'}$. If $\nu,\nu'\in \Dr(M)$, then $\underline{\nu}=\mathcal{B}=\underline{\nu'}$, so that then 
$$\nu\sim\nu'~\Longleftrightarrow[\overline{\nu}]=[\overline{\nu'}]$$
The set of bases $\mathcal{B}$ of the underlying matroid $M$ restricts $[\overline{\nu}]$, as follows.
\begin{lemma} Let $M=(E,\mathcal{B})$ be a matroid, let $\nu\in \Dr(M)$, and let $(S,ab,cd)\in Z(r,E)$. If $Sab\not\in \mathcal{B}$, then $(S,ab,cd)\in [\overline{\nu}]$.
\end{lemma}
\proof If $Sab\not\in \mathcal{B}$, then $\overline{\nu}(Sab)+\overline{\nu}(Scd)=\infty$. Then the minimum of 
$$\overline{\nu}(Sab)+\overline{\nu}(Scd), ~ \overline{\nu}(Sac)+\overline{\nu}(Sbd), ~ \overline{\nu}(Sad)+\overline{\nu}(Sbc)$$
equals $\infty$,  or is  attained by the latter two expressions. Then $\nu(Sac)+\nu(Sbd)=\nu(Sad)+\nu(Sbc)$ either way.\endproof
We introduce notation that facilitates exploiting this lemma. For any matroid $M=(E,\mathcal{B})$, let
\begin{enumerate}
\item $Z(M):=\{(S,ab,cd)\in Z(r,E):  Sac, Sbd, Sad, Sbc\in \mathcal{B}\}$,
\item $Z_0(M):=\{(S,ab,cd)\in Z(M): Sab\not\in \mathcal{B}\text{ or } Scd\not\in \mathcal{B}\}$,
\item $Z_1(M):= Z(M)\setminus Z_0(M)$.
\end{enumerate}
and for any $\nu\in \Dr(M)$, put
$$[\nu]:=[\overline{\nu}]\cap Z(M), \qquad [\nu]_1:=[\overline{\nu}]\cap Z_1(M)$$
Note that this does not clash with our earlier definition of $[\nu]$ for $\nu:\binom{E}{r}\rightarrow \Rinf$: in the only possibly conflicting case that $\nu:\binom{E}{r}\rightarrow \R$, $M$ must be a uniform matroid and then $Z(M)=Z(r,E)$.
\begin{corollary}\label{cor:eq} For any $\nu\in \Dr(M)$, we have $Z(r,E)\setminus [\overline\nu]\subseteq Z_1(M)$. Hence 
$$\nu\sim\nu' \Longleftrightarrow [\nu]=[\nu'] \Longleftrightarrow [\nu]_1=[\nu']_1$$
\end{corollary}

The Dressian $\Dr(M)$ is the support of a fan
$$\dr(M):=\left\{\overline{D(\nu)}: \nu\in \Dr(M)\right\}$$
whose cells are the topological closures of combinatorial equivalence classes $$D(\nu):=\{\nu'\in \R^{\mathcal{B}}: \nu'\sim \nu\}.$$ Technically speaking the fan $\dr(M)$ and the Dressian $\Dr(M)\subseteq \R^{\mathcal{B}}$ are different objects, but we will nevertheless refer to the elements of the fan  $\dr(M)$ as  cells of the Dressian. 

Corollary \ref{cor:eq} gives a first upper bound on $\#\dr(M)$. Since for any $\nu, \nu'\in \Dr(M)$ we have
$$D(\nu)\neq D(\nu')\Longleftrightarrow \nu\not\sim\nu'\Longleftrightarrow  [\nu]_1\neq [\nu']_1,$$
so that $\#\dr(M)\leq \#\{[\nu]_1: \nu\in \Dr(M)\}\leq \#\{Z: Z\subseteq Z_1(M)\}$. Hence
$$\log_2 \#\dr(M)\leq |Z_1(M)|\leq |Z(r,E)|=\binom{n}{r-2}\cdot \binom{n-r+2}{4}\cdot 6\leq \binom{n}{r}\frac{(n-r)^2r^2}{4}$$
for each matroid $M$ of rank $r$ on $n$ elements.

\subsection{Residue matroids and matroid subdivisions}
Let $M=(E,\mathcal{B})$ be a matroid. For any $\nu\in \Dr(\nu)$, the {\em initial} or {\em residue} matroid of $\nu$ is the matroid $M_0(\nu):=(E, \mathcal{B}_0(\nu))$ with bases $$\mathcal{B}_0(\nu):=\left\{B\in \mathcal{B}: \nu(B)=\min_{B'\in \mathcal{B}} \nu(B')\right\}.$$
For any $\nu\in \Dr(M)$ and $w\in\R^E$, we define $\nu^w:\mathcal{B}\rightarrow \R$ by setting
$$\nu^w(B):=\nu(B)+\sum_{e\in B} w(e)$$
for all $B\in \mathcal{B}$. It is straightforward that $\nu^w\in \Dr(M)$ and $[\nu^w]=[\nu]$.
\begin{lemma} Let $M=(E,\mathcal{B})$ be a matroid and let $\nu\in \Dr(M)$. The following are equivalent for any $(S, ab,cd)\in Z_1(M)$:
\begin{enumerate}
\item $(S, ab,cd)\in [\nu]$; and
\item for each $w\in\R^E$, we have $Sac, Sbd\in \mathcal{B}_0(\nu^w)$ if and only if $Sad, Sbc\in \mathcal{B}_0(\nu^w)$.
\end{enumerate}
\end{lemma}
\ignore{
\proof Suppose $(S, ab,cd)\in [\nu]$. Let $w\in\R^E$ be such that $Sac, Sbd\in \mathcal{B}_0(\nu^w)$. Without loss of generality, we may assume $w=0$. Then 
$\nu(Sac)=\nu(Sbd)=c:=\min_{B'\in \mathcal{B}} \nu(B').$
As $$2c=\nu(Sac)+\nu(Sbd)=\nu(Sad)+\nu(Sbc)\geq 2c,$$
it follows that $\nu(Sad)=\nu(Sbc)=c=\min_{B'\in \mathcal{B}} \nu(B')$ as well, as required.

For the reverse implication, suppose  that $(S, ab,cd)\not\in [\nu]$. Then $(S, \{ad, bc\})\in [\nu]$, say, and we have
$$\nu(Sab)+\nu(Scd) = \nu(Sac)+\nu(Sbd) < \nu(Sad)+\nu(Sbc).$$
We will construct a $w\in\R^E$ so that $Sac, Sbd\in \mathcal{B}_0(\nu^w)$ but  $Sad, Sbc\not\in \mathcal{B}_0(\nu^w)$.
For a sufficiently large value of $W$, put $w(e)=-W$ for all $e\in S$ and $w(e)=W$ for all $e\in E\setminus Sabcd$, so that
 $$S\subseteq B\subseteq E\setminus (Sabcd)$$ for $B\in \mathcal{B}_0(\nu^w)$. 
Pick $w(e)\in \R$ for $e\in \{a,b,c,d\}$ so that $\nu^w(Sac)=\nu^w(Sbd)=\nu^w(Sad)=\nu^w(Sbc)$. As $(S, ab,cd)\not\in [\nu^w]=[\nu]$ and $(S, \{ad, bc\})\in [\nu^w]=[\nu]$, it follows that 
$$\nu^w(Sab)+\nu^w(Scd) = \nu^w(Sac)+\nu^w(Sbd) < \nu^w(Sad)+\nu^w(Sbc).$$
Then $Sac, Sbd\in \mathcal{B}_0(\nu^w)$ but  $Sad, Sbc\not\in \mathcal{B}_0(\nu^w)$, as required.
\endproof
}
The {\em matroid polytope} of a matroid $M$ on $E$ is 
$$P_M:=\mbox{conv.hull}\left\{e_B: B\in \mathcal{B}\right\}\subseteq \R^E$$
where $e_B:=\sum_{b\in B} e_b\in\{0,1\}^E$ is the incidence vector of the set $B\subseteq E$. If $\nu\in \Dr(M)$ then the set
$$\mathcal{P}(\nu):=\left\{P_{M_0(\nu^w)}:w\in \R^E\right\}.$$
is a polyhedral complex whose cells are matroid polyhedra (see \cite[Proposition 2.2]{Speyer2008}). Combinatorial equivalence of valuations amounts to equality of the corresponding polyhedral complexes.
\begin{lemma} Let $M$ be a matroid and let $\nu, \nu'\in \Dr(M)$. Then 
$$\mathcal{P}(\nu)=\mathcal{P}(\nu') ~\Longleftrightarrow ~[\nu]=[\nu'].$$
\end{lemma}
\ignore{
\proof Necessity:  for each $(S, ab,cd)\in Z_1(M)$, condition (2) of Lemma {lem:residue}  depends only on $\mathcal{P}(\nu)$. Hence if $\mathcal{P}(\nu)=\mathcal{P}(\nu')$, then $[\nu]=[\nu']$.

Sufficiency:  As each $M_0(\nu^w)$ is determined by $[\nu]$ and $w$ alone, $\mathcal{P}(\nu)$ is determined by only $[\nu]$. \endproof
}
In other words, each cell of the Dressian $\dr(M)$ corresponds to a unique matroid decomposition of the matroid polytope of $M$

The {\em spread} of a polyhedral complex $\mathcal{P}$, $\mbox{\em spread}(\mathcal{P})$, is defined as  the number of inclusion-wise maximal cells of $\mathcal{P}$. Speyer showed \cite{Speyer2008}:
\begin{theorem} Let $M$ be a matroid of rank $r$ on $n$ elements, and let $\nu\in \Dr(M)$. Then $$\mbox{\em spread}(\mathcal{P}(\nu))\leq \binom{n-2}{r-2}.$$
\end{theorem}
In \cite{HJJS2009}, it is noted that the dimension of the Dressian is bounded in terms of the spread.
\begin{corollary}\label{cor:spreaddim}Let $M$ be a matroid of rank $r$ on $n$ elements. Then 
$\dim \dr(M)\leq \binom{n-2}{r-2}+n-1$.
\end{corollary}
\proof For each $\nu\in \Dr(M)$, we have $\dim  D(\nu)\leq \mbox{\em spread}(\mathcal{P}(\nu))+n-1$.\endproof

\subsection{Valuated matroids of rank 2}
In \cite[Section 4]{SpeyerSturmfels2004}, Speyer and Sturmfels identify the Dressian in rank 2 with the space of phylogenetic trees. Slightly extending to the case when the underlying matroid $M$ is not uniform, their results can be paraphrased as follows. 
\begin{theorem} Let $M=(E, \mathcal{B})$ be a matroid of rank 2. The following are equivalent for a function $\nu:\mathcal{B}\rightarrow \R$: 
\begin{enumerate}
\item $\nu\in \Dr(M)$.
\item There is a tree $T$ with leafs $E\subseteq V(T)$ and a function $\ell: E(T)\rightarrow \R$ so that
\begin{itemize}
\item $\ell(e)>0$ for each internal edge $e$ of $T$;
\item if $\{a,b\}\in \binom{E}{2}\setminus\mathcal{B}$ then $a$ and $b$ have a common neighbor in $T$; and
\item if $\{a,b\}\in \mathcal{B}$, then $\nu(ab)=\ell[E(P)]$ where $P$ is the unique path from $a$ to $b$ in $T$.
\end{itemize}
\end{enumerate}
\end{theorem}
Here $V(T)$ is the set of vertices of $T$ and $E(T), E(P)$ denote the set of edges of $T$ and $P$ respectively, and $\ell[E(P)]:=\sum_{e\in E(P)} \ell(e)$ is the total $\ell$-length of the path $P$. 

It is not difficult to see that if a valuation $\nu$ corresponds to a metric tree $(T,\ell)$ as in the theorem, then the combinatorial type $[\nu]$ is fully determined by the topology of $T$, and that if $\nu'$ corresponds to $(T', \ell')$ then $\nu \sim \nu'$  if and only if $T$ and $T'$ are topologically equivalent. Thus cells of the Dressian of a matroid of rank 2 on $E$ correspond 1-1 to topological trees with set of leafs $E$; the dimension of each cell is bounded by the number of arcs of the corresponding tree, and the number of cells equals the number of topological trees with leafs $E$. 

A somewhat more detailed but straightforward analysis yields:
\begin{corollary}\label{cor:treebounds} Let $M$ be a matroid of rank 2 on $n$ elements, with $t$ parallel classes. Then
$$\dim \dr(M)\leq n+t-3\text{ and } \#\dr(M)\leq 2^t t^n.$$
\end{corollary}

\section{The dimension of the Dressian}
If $X$ is a finite set, then an {\em exact $k$-cover of $X$} is a multiset $\mathcal{A}=\{A_1,\ldots, A_m\}$ of nonempty subsets $A_i\subseteq X$ so that  $\#\{i: x\in A_i\}=k$ for each $x\in X$. If $L\subseteq \mathbb{F}^X$, then we write $$L_A:=\{(x_a: a\in A): x \in L\}$$ for the projection of $L$ on $\mathbb{F}^A$.
\begin{lemma} \label{lem:notshearer}Let $\mathbb{F}$ be a field, let $X$ be a finite set, and let  $\mathcal{A}$ be an exact $k$-cover of $X$. 
Then for each linear subspace  $L\subseteq \mathbb{F}^X$, we have
$$\dim(L)\leq \frac{1}{k} \sum_{A\in \mathcal{A}} \dim(L_A).$$
\end{lemma}
\proof For any $A\subseteq X$, denote $L^A=\{x\in L: x_a=0\text{ for all } a\in A\}$. Then 
for any two sets $A, A'\subseteq X$, we have $L^A\cap L^{A'} = L^{A\cup A'}$ and $L^A+ L^{A'}\subseteq L^{A\cap A'}$. Hence
$$\dim\left(L^A\right)+\dim\left(L^{A'}\right)=\dim\left(L^A\cap L^{A'} \right)+\dim\left(L^A + L^{A'}\right)\leq  \dim\left(L^{A\cup A'}\right)+\dim\left( L^{A\cap A'}\right)$$
As $\dim(L) = \dim(L_A)+\dim(L^A)$ we obtain
$$\dim(L_A)+\dim(L_{A'})\geq  \dim(L_{A\cup A'})+\dim( L_{A\cap A'})$$
Without loss of generality, we may assume that  $\mathcal{A}$ is one of the exact $k$-cover of $X$ such that  $\sum_{A\in \mathcal{A}} \dim(L_A)$ is as small as possible. Among these minimal exact $k$-covers, choose $\mathcal{A}$ so that  $\sum_{A\in \mathcal{A}} |A|^2$ is as large as possible. 

We claim that $\mathcal{A}$ consists of $k$ multiples of $X$. If not, let $A\in \mathcal{A}$, $A\neq X$. Considering any $x\in X\setminus A$, there must be an $A'\in \mathcal{A}$ so that $x\in A'\not\supseteq A$, as  $\mathcal{A}$ is an exact $k$-cover. So $A\setminus A'$ and $A'\setminus A$ are nonempty. 
Let $\mathcal{A}'$ arise from $\mathcal{A}$ by replacing $A, A'\in \mathcal{A}$ with $A\cap A'$, $A\cup A'$. Then 
$$\left(\sum_{A\in \mathcal{A}'} \dim(L_A)\right) - \left(\sum_{A\in \mathcal{A}} \dim(L_A)\right)=\left(\dim(L_{A\cup A'})+\dim( L_{A\cap A'})\right)-\left(\dim(L_A)+\dim(L_{A'})\right)\leq 0$$
and 
$$\left(\sum_{A\in \mathcal{A}'} |A|^2\right)-\left(\sum_{A\in \mathcal{A}} |A|^2\right)=\left(|A\cap A'|^2+|A\cup A'|^2\right)-\left(|A|^2+|A'|^2\right)>0$$
as $A\setminus A'$ and $A'\setminus A$ are nonempty. This contradiction proves the claim. 

It follows that $\sum_{A\in \mathcal{A}} \dim(L_A)=k\dim(L)$, which proves the lemma.\endproof

We use Lemma \ref{lem:notshearer} to prove upper bounds on the dimension of the Dressian of rank-$r$ matroids based on similar bounds in lower rank. 
\begin{lemma} \label{lem:proj} Let $M=(E,\mathcal{B})$ be a matroid of rank $r$, and let $S\subseteq E$ be an independent set of $M$. For each $D\in \dr(M)$ there is a $D'\in \dr(M/S)$ so that $D/S\subseteq D'$.
\end{lemma}
\proof It suffices to show that any two elements of $D/S$ are combinatorially equivalent, that is, are contained in the same equivalence class $D'$. So let $\nu, \nu'\in D$, and consider $\nu/S, \nu'/S\in D/S$. If $\nu/S\not\sim \nu'/S$, then there is a symbol $(S', ab,cd)\in [\nu/S]\setminus [\nu'/S]$, say, and then trivially $(S\cup S', ab,cd)\in [\nu]\setminus [\nu']$. But this would contradict that $\nu, \nu'\in D\in \dr(M)$ and that hence $[\nu]= [\nu']$.\endproof

\begin{theorem} Let $M=(E,\mathcal{B})$ be a matroid of rank $r$, on $n:=|E|$ elements, and let $t\leq r$. Then
 \begin{equation}\label{eq:notentropy}
 \frac{\dim \dr(M)}{\binom{n}{r}}\leq \frac{\max\{\dim \dr(M/S): ~S\in \binom{E}{r-t}\text{ independent in }M\}}{\binom{n-r+t}{t}}.
 \end{equation}
 \end{theorem}
\proof Put $s:=r-t$. Let $\mathcal{S}^*:=\left\{S\in \binom{E}{s}: S\text{ independent in }M\right\}$ and for each $S\in\mathcal{S}^*$, let 
$$A(S):=\{B\in \mathcal{B}: S\subseteq B\}.$$ Then $\mathcal{A}=\left\{A(S): S\in \mathcal{S}^* \right\}$ is an exact $\binom{r}{s}$-cover of $\mathcal{B}$. 

Let $D\in \Dr(M)$. For each $S\in \mathcal{S}^*$, we have 
$$D_{A(S)}=\{(\nu(B): B\supseteq S): \nu\in D\}\cong \{\nu/S: \nu\in D\}=:D/S$$
By Lemma \ref{lem:proj}, there is a $D'\in \dr(M/S)$ so that $D/S\subseteq D'$.
Taking $L=\mbox{lin.hull}(D)$, it follows that 
$L_{A(S)}\subseteq \mbox{lin.hull}(D_{A(S)})$ and in particular, 
$$\dim(L_{A(S)})\leq \dim(D/S)\leq \dim(D')\leq  \binom{n-s}{r-s}\alpha$$
where $\alpha$ equals the RHS of \eqref{eq:notentropy}. 
By Lemma \ref{lem:notshearer}, we have
$$\dim(D)=\dim(L)\leq \frac{1}{\binom{r}{s}}\sum_{S\in \mathcal{S}^*} \dim(L_{A(S)})\leq  \frac{\binom{n}{s}}{\binom{r}{s}} \alpha \binom{n-s}{r-s}=\binom{n}{r}\alpha.$$
The theorem follows.
\endproof

\begin{corollary} \label{cor:dim}Let $M=(E,\mathcal{B})$ be a matroid of rank $r\geq 3$ on $n$ elements. Then
$$\dim \dr(M)\leq \binom{n}{r}\frac{3}{n-r+3}$$
\end{corollary}
\proof For each  $S\in \binom{E}{n-r+3}$ so that $S$ is independent in $M$, the matroid $N=M/S$ has rank 3 and $k=n-r+3$ elements, so that
$$\dim \dr(M/S)\leq \binom{k-2}{2}+k-1= \binom{n-r+2}{2}$$ by Corollary \ref{cor:spreaddim}.  Then
$$\frac{\dim\dr(M)}{\binom{n-r+3}{3}}\leq  \frac{\max\{\dim \dr(M/S): ~S\in \binom{E}{r-t}\text{ independent in }M\}}{\binom{n-r+3}{3}}\leq \frac{\binom{n-r+2}{2}}{\binom{n-r+3}{3}}\leq \frac{3}{n-r+3}$$
by applying the theorem with $t=3$.\endproof
Note that applying the above argument with $t=2$ or $t\geq 4$ would yield a strictly weaker upper bound.

\section{The number of  cells of the Dressian}

\subsection{Encoding cells as linear subspaces} 
Let $M=(E,\mathcal{B})$ be a matroid. For each $Z\subseteq Z(M)$, we will write
 $$L(Z):=\left\{\nu \in \R^{\binom{E}{r}}: \nu(Sac)+\nu(Sbd)=\nu(Sad)+\nu(Sbc) \text{ for }(S,ab,cd)\in Z, \nu(B)=0 \text{ for } B\not\in \mathcal{B}\right\}.$$
For any matroid $M$, we define $U(M):=L(Z_0(M))$.

The linear space $L([\nu])$ uniquely determines the combinatorial class $[\nu]$. 
\begin{lemma}\label{lem:bijection} Let $M=(E, \mathcal{B})$ be a matroid. If $\nu, \nu'\in \Dr(M)$ then 
$$D(\nu)=D(\nu')\Longleftrightarrow \nu\sim \nu'  \Longleftrightarrow L([\nu])=L([\nu'])$$
\end{lemma}
\proof By definition $D(\nu)=D(\nu')\Longleftrightarrow \nu\sim \nu'    \Longrightarrow L([\nu])=L([\nu'])$. If  $\nu\not\sim \nu'$, then there exists a symbol $(S, ab, cd)\in [\nu]\setminus [\nu']$ say. Then $\nu''\in L([\nu'])\setminus L([\nu])$, where $\nu''\in \R^{\binom{E}{r}}$ arises by extending $\nu'\in \R^{\mathcal{B}}$ with zeros, since
$$\nu''(Sac)+\nu''(Sbd)=\nu'(Sac)+\nu'(Sbd)\neq\nu'(Sad)+\nu'(Sbc)=\nu''(Sad)+\nu''(Sbc).$$
 It follows that $L([\nu])=L([\nu'])\Longrightarrow  \nu\sim \nu' $.
\endproof
It follows that  $\#\dr(M)=\#\mathcal{L}(M)$ where $$\mathcal{L}(M):=\{L([\nu]): \nu \in\Dr(M)\}.$$ For technical reasons, we also define
$\mathcal{L}(M, U):=\{L\in \mathcal{L}(M): L\subseteq U\}.$
The fact that each $L\in \mathcal{L}(M, U)$ is a linear subspace of $U$ determined by linear equations from a bounded set yields the following upper bound on the number of elements of $\#\mathcal{L}(M, U)$. 
\begin{lemma} \label{lem:dimbound} Let $M=(E,\mathcal{B})$ be a matroid, and let $U\subseteq \R^{\binom{E}{r}}$ be a linear subspace of dimension $u$. Then 
$$\#\mathcal{L}(M,U)\leq \binom{|Z_1(M)|}{\leq u}. $$
\end{lemma}
\proof We may assume that $U\subseteq U(M)$. Then
\begin{eqnarray*}\mathcal{L}(M,U)=& \{L([\nu]): \nu\in \Dr(M), L([\nu])\subseteq U\}&\subseteq \\
\subseteq&\{L(C): Z_0(M)\subseteq C\subseteq Z(M), L(C)\subseteq U\}&\subseteq \\
\subseteq&\{L(C)\cap U: Z_0(M)\subseteq C\subseteq Z(M)\}&=\\
=&\{L(C_1) \cap U: C_1\subseteq Z_1(M)\}&=\\
=&\{L(C_1)\cap U: C_1\subseteq Z_1(M), |C_1|\leq u\}
\end{eqnarray*}
The last equality follows considering that
$$L(C_1)\cap U=\{\nu \in U: \nu(Sac)+\nu(Sbd)=\nu(Sad)+\nu(Sbc) \text{ for all }(S,ab,cd)\in C_1\}.$$ 
Then for each $C_1\subseteq Z_1(M)$ there is a $C_1'\subseteq C_1$ with $|C_1'|\leq u$ so that $L(C_1)\cap U=L(C_1')\cap U$.
It follows that $$\#\mathcal{L}(M,U)\leq \#\{C_1: C_1\subseteq Z_1(M), |C_1|\leq u\}\leq \binom{|Z_1(M)|}{\leq u}$$
as required.
\endproof

\subsection{Subspace containers}
Let  $\mathcal{L},\mathcal{C}$ be collections of linear subspaces from $\mathbb{F}^X$, and let $\alpha, \beta\geq 0$. 
We say that   $\mathcal{C}$ is an {\em $(\alpha, \beta)$-container} for $\mathcal{L}$ if 
\begin{enumerate}
\item for each $L\in \mathcal{L}$ there is an $C\in \mathcal{C}$ so that $L\subseteq C$, and 
\item $\dim(C)\leq \alpha |X|$ for all $C\in \mathcal{C}$, and
\item $\ln \#\mathcal{C}\leq \beta |X|$
\end{enumerate}
A container for $\mathcal{L}(M)$ gives a bound on $\#\dr(M)=\#\mathcal{L}(M)$.
\begin{lemma}\label{lem:containercount} Let $M$ be a matroid or rank $r$ on $n$ elements. If $\mathcal{L}(M)$ has an $(\alpha, \beta)$-container, then
$$\ln\#\dr(M)\leq \binom{n}{r}\left(\alpha\ln\left(\frac{n^4}{\alpha}\right)+\beta\right).$$
\end{lemma}
\proof Let $\mathcal{C}$ be an $(\alpha, \beta)$-container for $\mathcal{L}(M)$. If $C\in \mathcal{C}$, then  $\dim(C)\leq \alpha\binom{n}{r}$ and hence
$$\#\mathcal{L}(M,C)\leq \binom{|Z_1(M)|}{\leq \dim(C)}\leq \binom{\binom{n}{r}n^4/4}{\leq \alpha\binom{n}{r}}$$
by  Lemma \ref{lem:dimbound}.  As $\mathcal{C}$ is a container for $\mathcal{L}(M)$ we have 
$\#\mathcal{L}(M)= \sum_{C\in \mathcal{C}}\#\mathcal{L}(M,C)$, hence
$$\ln \#\mathcal{L}(M)\leq \ln\#\mathcal{C} +\ln\max_{C\in \mathcal{C}}\#\mathcal{L}(M,C)\leq \beta\binom{n}{r}+\alpha\binom{n}{r}\ln\left(\frac{n^4}{\alpha}\right).$$
By Lemma \ref{lem:bijection}, we have $\ln\#\dr(M)=\ln \#\mathcal{L}(M)$, and the current lemma follows.
\endproof
\begin{theorem} \label{thm:subspacebound} Let $M$ be a matroid or rank $r$ on $n$ elements, and let $u:=\dim(U(M))$. Then 
$$\ln \#\dr(M)\leq u\ln\left(\binom{n}{r}n^4/u\right).$$
\end{theorem}
\proof For each $\nu\in\Dr(M)$, we have $[\nu]\supseteq Z_0(M)$ and hence $L([\nu])\subseteq L(Z_0(M))=:U(M)$. Then $\mathcal{C}=\{U(M)\}$ is evidently an $(\alpha, \beta)$ container for $\mathcal{L}(M)$, with
$$\alpha:=\frac{u}{\binom{n}{r}}, \qquad \beta=0.$$
By Lemma \ref{lem:containercount}, it follows that
$\ln\#\dr(M)\leq \binom{n}{r}\left(\alpha\ln\left(\frac{n^4}{\alpha}\right)+\beta\right)= u\ln\left(\binom{n}{r}n^4/u\right)$
as required.
\endproof
If $\mathcal{L}$ is a set of linear subspaces of $\mathbb{F}^X$ and  $A\subseteq X$, we write $\mathcal{L}_A:=\{L_A: L\in \mathcal{L}\}$.
\begin{lemma} \label{lem:minor}Let $M=(E, \mathcal{B}$ be a matroid of rank $r$ on $n$ elements, let $S\subseteq E$ be a set of cardinality $s\leq r$, and let $A:=\left\{B\in \binom{E}{r}: B\supseteq S\right\}$. Then $\mathcal{L}(M)_A$ has an $(\alpha, \beta)$-container, where 
$$\alpha:=\frac{\dim\dr(M/S)}{\binom{n-s}{r-s}}, \qquad \beta:=\frac{\ln\#\dr(M/S)}{\binom{n-s}{r-s}}.$$
\end{lemma}
\proof If $S$ is a dependent set of $M$, then $\mathcal{L}(M)_A$ only contains the trivial subspace $\{0\}$, and hence has a $(0,0)$-container. So $S$ is independent.
We claim that $\mathcal{L}(M/S)$ is an $(\alpha,\beta)$-container for $\mathcal{L}(M)/S$. As $\mathcal{L}(M)/S$ is a set of subspaces from $\R^{\binom{E\setminus S}{r-s}}$, we need to argue that 
\begin{enumerate}
\item for each $L\in \mathcal{L}(M)/S$ there is an $C\in \mathcal{L}(M/S)$ so that $L\subseteq C$, and 
\item $\dim(C)\leq \alpha \binom{n-s}{r-s}$ for all $C\in \mathcal{L}(M/S)$, and
\item $\ln \#\mathcal{L}(M/S)\leq \beta \binom{n-s}{r-s}$.
\end{enumerate}
To see (1), note that each $L\in \mathcal{L}(M)$ is of the form $L=L([\nu])$ for some $\nu\in \Dr(M)$, so that
$$L/S=L([\nu])/S\subseteq L([\nu/S])=:C\in \mathcal{L}(M/S)$$
as required.
To show (2), we will argue that for each $C\in \mathcal{L}(M/S)$ there is a $D\in \dr(M/S)$ so that 
$\dim(C)\leq\dim(D)\leq \dim \dr(M/S)$. To see this, let $\nu \in \Dr(M/S)$ be chosen so that $L=L([\nu])$, and let $D=D(\nu)$. Consider a vector $\mu\in L([\nu])_{\mathcal{B}}$. Then $[\nu+\epsilon\mu]=[\nu]$ for all sufficiently small $\epsilon>0$, so that $\nu, \nu+\epsilon\mu\in D$, and hence $\mu\in \mbox{lin.hull}(D)$. Then
$$\dim(C)=\dim (L([\nu]))=\dim (L([\nu])_{\mathcal{B}})\leq \dim (D)$$
as required.
Finally, (3) is straightforward as $ \#\mathcal{L}(M/S)=\#\dr(M/S)$.

As $\mathcal{L}(M)_{A(S)}$ is isomorphic to $\mathcal{L}(M)/S$ via the coordinate map $A(S)\rightarrow \binom{E-S}{r-s}: B\mapsto B\setminus S$, it follows that $\mathcal{L}(M)_{A(S)}$ has an $(\alpha, \beta)$-container as well.
\endproof

\begin{lemma} \label{lem:container}Let $\mathbb{F}$ be a field, let $X$ be a finite set, let  $\mathcal{A}$ be an exact $k$-cover of $X$, and let $\alpha, \beta\geq 0$. 
Let  $U\subseteq\mathbb{F}^X$ be a linear subspace such that $\dim(U)\geq 2\alpha |X|$, and let $\mathcal{L}$  be a collection of linear subspaces of $U$. If $\mathcal{L}_A$ has an $(\alpha,\beta)$-container for each $A\in \mathcal{A}$, then $\mathcal{L}$ has an $(\alpha', \beta')$-container, where
$$\alpha':=2\alpha+\frac{\max\{|A|: A\in \mathcal{A}\}}{|X|} \text{ and } \beta':=\beta\ln\left(\frac{\dim(U)}{\alpha|X|}-1\right)$$
\end{lemma}

\proof We prove the lemma by induction on $u=\dim(U)$. If $2\alpha|X|\leq \dim(U)\leq \alpha'|X|$, then $\mathcal{C}=\{U\}$ is an $(\alpha', 0)$-container for $\mathcal{L}$. If $\dim(U)>\alpha'|X|$, let $A^*\in \mathcal{A}$ be such that 
$$\frac{\dim U_{A^*}}{|A^*|}=\alpha^*:=\max\left\{\frac{\dim (U_A)}{|A|}: A\in \mathcal{A}\right\}$$
Let $\mathcal{C}^*=\{C_1,\ldots ,C_k\}$ be an $(\alpha,\beta)$-container for $\mathcal{L}_{A^*}$. 
For $i=1,\ldots, k$, let $$\mathcal{L}^i:=\{L\in \mathcal{L}: L_A\subseteq C_i\}.$$ 
Then $\mathcal{L}^i$ is a collection of subspaces from  $U^i:=\{u\in U: u_A\in C_i\}$, where 
$$\dim(U^i)\geq \dim(U)-|A|\geq \alpha'|X|-\max\{|A|: A\in \mathcal{A}\}=2\alpha|X|.$$
By induction, there is a $(\alpha',\beta^i)$-container $\mathcal{C}^i$ for $\mathcal{L}^i$ such that $\beta^i= \beta f(\dim(U_i))$, where
$$ f(x):= \ln\left(\frac{x}{\alpha|X|}-1\right).$$
We will show that $\mathcal{C}:=\mathcal{C}^1\cup \cdots\cup \mathcal{C}^k$ is an $(\alpha',\beta')$-container for $\mathcal{L}$, that is, 
\begin{enumerate}
\item for each $L\in \mathcal{L}$ there is an $C\in \mathcal{C}$ so that $L\subseteq C$, and 
\item $\dim(C)\leq \alpha' |X|$ for all $C\in \mathcal{C}$, and
\item $\ln \#\mathcal{C}\leq \beta' |X|$
\end{enumerate}
To see (1), note that for each $L\in \mathcal{L}$ there is a $i$ so that $L\in \mathcal{L}^i$, and then $L\subseteq C$ for a $C\in \mathcal{C}^i\subseteq \mathcal{C}$, as required. 
Condition (2) is also true by construction, since if $C\in \mathcal{C}$, then there is a $i$ so that $C\in \mathcal{C}^i$, and then $\dim(C)\leq \alpha' |X|$. 
It remains to show (3). Using Lemma \ref{lem:notshearer}, we find that
$$\dim(U)\leq \frac{1}{k}\sum_{A\in\mathcal{A}} \dim(U_A)\leq \frac{1}{k}\sum_{A\in\mathcal{A}} \alpha^*|A|=\alpha^*|X|.$$
Then $\dim(U)/|X|\leq \alpha^*=\dim(U_{A^*})/|A^*|$. It follows that
$$\dim(U)-\dim(U^i)= \dim(U_{A^*})-\dim((U^i)_{A^*})\geq (\alpha^*-\alpha)|A^*|\geq\left(\frac{\dim(U)}{|X|}-\alpha\right) |A^*|=:t$$
Then for each $i$, we have $\dim(U^i)\leq u-t$ and hence 
$$f(\dim(U^i))\leq f(u-t)\leq f(u)-f'(u)t= f(u)-\frac{t}{u-\alpha |X|}=f(u)-\frac{|A^*|}{|X|}$$ 
as $f$ is an increasing and concave function. By construction of $\mathcal{C}$, we have
$$\ln \#\mathcal{C}\leq \ln \#\mathcal{C}_{A^*}+\max_i\ln \#\mathcal{C}^i\leq \beta |A^*|+ \beta f(u')|X|\leq \beta f(u)|X|=\beta'|X|,$$
as required.
\endproof
Applying this lemma to the coordinate projections of Lemma \ref{lem:minor}, we obtain the following technical statement.
\begin{lemma}\label{lem:monty}Let $M=(E, \mathcal{B})$ be a matroid of rank $r\geq 2$ on $n:=|E|$ elements, and  let $2\leq s \leq r$. 
Let 
$$\alpha^*:= \frac{\max\{\dim\dr(M/S): S\in \binom{E}{s}\}}{\binom{n-s}{r-s}}, \qquad\beta^*:=\frac{\max\{\ln\#\dr(M/S):S\in \binom{E}{s}\}}{\binom{n-s}{r-s}}.$$
Then
$\mathcal{L}(M)$ has an $(\alpha, \beta)$-container, 
where $$\alpha:=2\alpha^*+\frac{\binom{n-s}{r-s}}{\binom{n}{r}} \text{ and } \beta:=\beta^*\ln\left(\frac{\upsilon}{\alpha^*}-1\right)$$
with $\upsilon:=\dim(U(M))/\binom{n}{r}$.
\end{lemma}
\proof For each $S\in \binom{E}{s}$, let $A(S):=\left\{B\in \binom{E}{r}: S\subseteq B\right\}$. Then $\mathcal{A}=\left\{A(S):  S\in \binom{E}{s}\right\}$ is an exact $\binom{r}{s}$-cover of $\binom{E}{r}$. By Lemma \ref{lem:minor}, for each $S\in \binom{E}{s}$ the projection  $\mathcal{L}(M)_{A(S)}$ has an $(\alpha^*, \beta^*)$-container.
By Lemma \ref{lem:container}, $\mathcal{L}(M)$ has an $(\alpha, \beta)$ container, as required.
\endproof
If we apply the lemma without further assumptions on the matroid $M$ or its minors,  we obtain the following general bound.
\begin{theorem}\label{thm:count} Let $M$ be a matroid of rank $r\geq 3$ on $n$ elements. Then 
$$\ln\#\dr(M)\leq \frac{55\ln(n)+4\ln(n)^2}{n}\binom{n}{r}$$
\end{theorem}
\proof  Let $s:=r-2$. Using Lemma \ref{cor:treebounds}, we find that 
$$\alpha^*:= \frac{\max\{\dim\dr(M/S): S\in \binom{E}{s}\}}{\binom{n-s}{r-s}}\leq \frac{2(n-r+2)-3}{\binom{n-r+2}{2}}\leq \frac{4}{n-r+2}$$
and 
$$ \beta^*:=\frac{\max\{\ln\#\dr(M/S):S\in \binom{E}{s}\}}{\binom{n-s}{r-s}}\leq \frac{(n-r+2)\ln(2(n-r+2))}{\binom{n-r+2}{2}}\leq \frac{2\ln(2(n-r+2))}{n-r+1}$$
Applying Lemma \ref{lem:monty}, we find that $\mathcal{L}(M)$ has an $(\alpha, \beta)$-container with  
$$\alpha\leq \frac{8}{n-r+2}+\frac{\binom{n-r+2}{2}}{\binom{n}{r}}\leq \frac{11}{n}, \qquad \beta\leq \frac{2\ln(2(n-r+2))}{n-r+1}\ln\left(\frac{n-r+2}{8}\right)\leq \frac{4\ln(n)^2}{n}.$$
From Lemma \ref{lem:containercount}, we then obtain
$$\ln\#\dr(M)/\binom{n}{r}\leq \alpha\ln\left(\frac{n^4}{\alpha}\right)+\beta\leq \frac{11}{n}\ln(n\cdot n^4)+\frac{4\ln(n)^2}{n}.$$
The theorem follows.
\endproof

\section{A construction, and some conjectures}
\subsection{Valuations from matroids}
Let $N$ be a matroid of rank $r$ on $E$, and let $r_N:2^E\rightarrow \N$ be the rank function of $N$. 
We define the function $\nu_N:\binom{E}{r}\rightarrow \R$ by 
$$\nu_N(X):=r-r_N(X)$$
for all $X\in \binom{E}{r}$. The following is well-known, and straightforward to verify. 
\begin{lemma} Let $N$ be a matroid of rank $r$ on $E$. Then $\nu_N$ is a valuation of $U(r,E)$.
\end{lemma}
Here $U(r,E)$ denotes the uniform matroid of rank $r$ on ground set $E$.

\subsection{Dimension of cells}
We consider that or any matroid $N$, the dimension of $D(\nu_N)$ as a lower bound on the dimension of the Dressian containing this cell. 

Recall that the {\em Johnson graph} $J(r, E)$ is the graph with vertices $\binom{E}{r}$ so that $X,Y\in  \binom{E}{r}$ are adjacent if and only if $|X\setminus Y|=1$. When $N=(E,\mathcal{B})$ is a matroid of rank $r$, we write  $\mathcal{U}(N):=\binom{E}{r}\setminus \mathcal{B}$ for the set of dependent $r$-sets of $N$. We denote the number of connected components of a graph $G$ by $c(G)$. For any matroid $N=(E, \mathcal{B})$, let $c(N)$ be the number of connected components of $J(r,E)\setminus\mathcal{B}$.
\begin{lemma} \label{lem:compdim} Let $N$ be a matroid. Then $\dim(D(\nu_N))\geq c(N)$.
\end{lemma}
\proof Put $k=c(N)$, and let  $C_1,\ldots, C_k\subseteq \binom{E}{r}$ be the vertex sets of the connected components of $J(r,E)\setminus\mathcal{B}$. Then $$N_i:=\left(E, J(r,E)\setminus C_i\right)$$ is a matroid for $i=1,\ldots, k$.
The vectors $\nu_{N_1},\ldots, \nu_{N_k}\in \R^{\binom{E}{r}}$ are linearly independent as they have pairwise disjoint supports $\{X: \nu_{N_i}(X)>0\}=C_i$. Therefore, the open cone
$$K:=\left\{\sum_{i=1}^k \lambda_i\nu_{N_i}:\lambda_i>0\right\}$$
has dimension $\dim(K)=k$. We have $\nu_N=\sum_{i=1}^k \nu_{N_i}\in K$ and $[\kappa]=[\nu_N]$ for each $\kappa\in K$, so that 
$K\subseteq D(\nu)$ and hence $k=\dim(K)\leq \dim D(\nu_N)$.
\endproof
There exist sparse paving matroids $N=(E,\mathcal{B})$ of rank $r$ on $n$ elements so that $c(J(r,E)\setminus\mathcal{B})\geq \binom{n}{r}/n$. Immediately, we obtain:
\begin{theorem}$\frac{1}{n}\binom{n}{r}\leq \dim \dr(U(r,n))$.
\end{theorem}
For any matroid $M$, there is a sparse paving matroid $N$ so that $c(N)\geq c(M)$.  We find it difficult to imagine a cell of maximal  dimension that does not arise from a sparse paving matroid.
\begin{conjecture} For each fixed $r\geq 4$ there is an  $n_r\in \N$ so that for each $n\geq n_r$ there is a sparse paving matroid $N$ of rank $r$ on $n$ elements so that $\dim \dr(U(r,n)) = \dim D(\nu_N)$.
\end{conjecture}
If true, we would obtain that  $\dim\dr(U(r,n))\leq \alpha(J(r,n))\leq \binom{n}{r}/(n-r+1)$ whenever $n\geq n_r$, and it should be possible to improve Corollary \ref{cor:dim} by using a more accurate bound in fixed rank $t=4$ (or even $t=3$) for the base case of the argument. 

\subsection{Number of cells}
We next consider the consequences of the construction of valuations from matroids for the number of cells of the Dressian. 
\begin{lemma}Let $M_1, M_2$ be sparse paving matroids of rank $r$ on $E$, $|E|\geq 5$, with $2\leq r\leq |E|-2$. Then 
$$\nu_{M_1}\sim\nu_{M_2}~~\Leftrightarrow ~~M_1=M_2$$
\end{lemma}
\proof Suppose that $M_1, M_2$ contradict the theorem, chosen so that $n:=|E|$ is as small as possible. Since $n\geq 5$, we have $r>2$ or $n-r>2$. By matroid duality, we may assume that $n-r>2$. 

Suppose that $n>5$. Since $M_1$ and $M_2$ are sparse paving, no element $e\in E$ is a coloop of $M_1$ or $M_2$, so that $r(M_1\setminus e)=r=r(M_2\setminus e)$.  Hence for each $e\in E$ we have 
$$\nu_{M_1\setminus e}=\nu_{M_1}\setminus e\sim \nu_{M_2}\setminus e=\nu_{M_2\setminus e}$$
so that by minimality $M_1\setminus e=M_2\setminus e$ for each $e\in E$. Then $M_1=M_2$, contradicting our assumption.

So $n=5$, and hence $r=2$. Up to matroid isomorphism, each sparse paving matroid of rank 2 on 5 elements  $E=\{a,b,c,d,e\}$ is isomorphic to one of the following:
\begin{itemize}
\item $N_1=U(2,5)$; then $|[\nu_{N_1}]|=15$.
\item $N_2=(E, \mathcal{B}_2)$ where $\mathcal{B}_2=\binom{E}{2}\setminus\{ab\}$; then $|[\nu_{N_2}]|=9$.
\item $N_3=(E, \mathcal{B}_3)$ where $\mathcal{B}_3=\binom{E}{2}\setminus\{ab, cd\}$; then $|[\nu_{N_2}]|=5$.
\end{itemize}
So by counting the number of elements of $[\nu_{M_1}]=[\nu_{M_2}]$, we may determine an $i$ so that $N_i$ is isomorphic to both $M_1$ and $M_2$. For each $i$, a straightforward case check shows that if $M_1\cong N_i\cong M_2$ and  $[\nu_{M_1}]=[\nu_{M_2}]$, then $M_1=M_2$.
\endproof
Let $s(r,n)$ denote the number of sparse paving matroids of rank $r$ on $E=\{1,\ldots, n\}$. Since distinct sparse paving matroids $M$ give inequivalent valuations $\nu_M$, we immediately obtain:
\begin{theorem}For all $n\geq r\geq 2$, we have $s(n,r)\leq \#\dr(U(r,n))$.
\end{theorem}
As $\log_2 s(r,n)\geq \binom{n}{r}/n$ for any $r\leq n$ and $\log_2 s(r,n)\geq \binom{n}{r}\log_2(n)/n(1+o(1))$ as $n\rightarrow \infty$, we see that there is still a gap with  upper bound $\log_2 \#\dr(U(r,n))\leq c\binom{n}{r}\log_2(n)^2/n$ from Theorem \ref{thm:count}, even in fixed rank. We conjecture that the lower bound is asymptotically accurate. We conclude with a consideration about the structure of most valuations that may be helpful to improve the upper bounds.

We say that $\sigma\in \Dr(M)$ is a {\em spike} if $\#\{B:\sigma(B)\neq 0\} = 1$, and $\nu\in \Dr(M)$ is {\em smooth} if there are no $\nu', \sigma\in \Dr(M)$ so that $\nu=\nu'+\sigma$, and $\sigma$ is a spike. We expect that the equivalence classes of nonsmooth valuations predominate among cells of the Dressian $\Dr(M)$ unless $M$ is very far from being uniform. 
\begin{conjecture}For each fixed $r\geq 4$, we have
$\frac{\ln\#\{D(\nu): ~\nu\in \Dr(U(r,n)),~ \nu\text{ smooth}\}}{\ln \#\dr(U(r,n))} \rightarrow 0\text{ as }n\rightarrow \infty.$
\end{conjecture}
Since each valuation $\nu$ can be written $\nu=\nu'+\nu_N$ where $\nu'$ is smooth and $N$ is sparse paving, it would follow that 
$\frac{\ln \#\dr(U(r,n))}{\ln s(n,r)}\rightarrow 1\text{ as }n\rightarrow \infty$
for each fixed $r\geq 4$. To prove the conjecture, it may suffice to show an upper bound on the number of smooth cells in base rank $t=3$ or 4, and then apply the methods of this paper to count smooth cells. The results  will depend on the values of  
$\alpha^*$ and $\beta^*$ used in Lemma \ref{lem:monty}, so that we would need a good upper bound on the dimension of smooth cells as well.

\section{Final remarks}
\subsection{The subspace lemma} Lemma \ref{lem:notshearer} is structurally similar to {\em Shearer's Entropy Lemma} and its corollaries, {\em Shearer's Counting Lemma} and the {\em Loomis-Whitney inequality}. Lemma \ref{lem:notshearer} and Shearer's Lemma have the following common generalization.
\begin{lemma} \label{lem:polyshearer}Let $X$ be a finite set, let $f:2^X\rightarrow \R$ be a polymatroid, and let  $\mathcal{A}$ be an exact $k$-cover of $X$. Then 
$$f(X)\leq \frac{1}{k} \sum_{A\in \mathcal{A}} f(A).$$
\end{lemma}
Lemma \ref{lem:notshearer} follows by taking $f(A):=\dim(L_A)$, which $f$ is in fact the rank function of a matroid. Shearers Lemma for a random variable $Y=(Y_x: x\in X)$  follows by taking $f(A):=H(Y_A)$. The case of Lemma \ref{lem:notshearer} when $\mathbb{F}$ is a finite field is a corollary to Shearers Lemma, taking $Y$ the random variable which picks a vector uniformly from the linear space $L$ (which has exactly $|\mathbb{F}|^{\dim(L)}$ elements). The case $\mathbb{F}=\mathbb{Q}$ the also follows by reducing it to the finite field case with some additional argumentation, but we did not see a way to derive Lemma \ref{lem:notshearer} from Shearer's Lemma in full generality.

\subsection{The container lemma}
The $(\alpha,\beta)$-containers used for counting subspaces in this paper were inspired by the containers used for counting stable sets in graphs. Both types of containers serve the same purpose: subdividing the objects being counted (linear subspaces resp. stable sets)  into relatively few subclasses of bounded dimension/size. In the proof of Lemma \ref{lem:container}, Lemma \ref{lem:notshearer} is used to show that if a subspace $L$ has high dimension then some coordinate projection of $L$ has high dimension. This is analogous to the {\em supersaturation} driving the container method: that any subgraph spanned by many vertices has a vertex of high degree. 

While this analogy exists, we do acknowledge that Lemma \ref{lem:container} is not nearly subtle enough to generalize the container method for stable sets in graphs. If we encode a stable $S$ set of $G=(V,E)$ as a subspace 
$$L^S:= \R^S\times \{0\}^{V\setminus S}$$
and apply Lemma \ref{lem:container} to $\mathcal{L}^G:=\{L^S: S \text{ is a stable set of }G=(V,E)\}$ with projections to neighborhoods of vertices $\{v\}\cup N(v)$, say,
we obtain a container with $c\binom{n}{r}\ln(n)^2/n$ subspaces of dimension $O(\binom{n}{r}/n)$ in case of the Johnson graph $G=J(r,n)$. This relatively weak, as the regular container method produces a container with $o(\binom{n}{r}/n)$ subsets of size $O(\binom{n}{r}/n)$.

\bibliographystyle{alpha}
\bibliography{math}

\end{document}